\documentclass[11pt]{amsart}
\usepackage[T1]{fontenc}
\usepackage{amsmath,amssymb,amsthm}
\usepackage{booktabs}
\usepackage{tikz}
\usepackage[colorlinks=true,linkcolor=black,citecolor=black,urlcolor=black]{hyperref}

\theoremstyle{plain}
\newtheorem{theorem}{Theorem}[section]
\newtheorem{lemma}[theorem]{Lemma}
\newtheorem{proposition}[theorem]{Proposition}
\newtheorem{corollary}[theorem]{Corollary}
\theoremstyle{definition}
\newtheorem{definition}[theorem]{Definition}
\newtheorem{example}[theorem]{Example}

\newtheorem{problem}[theorem]{Problem}

\newcommand{\edgeless}[1]{\overline{K_{#1}}}
\newcommand{\ceil}[1]{\left\lceil #1\right\rceil}
\newcommand{\eg}{\operatorname{eg}}
\newcommand{\apex}{z}
\newcommand{\Pet}{\mathrm{Pe}}

\begin{document}

\title[Orientable genus of generalized corona graphs]{The orientable genus of generalized corona graphs}

\author{Jyoti Anant Pulgam$^1$}
\thanks{$^1$ Department of Mathematics and Statistics, Dr.\ Vishwanath Karad MIT World Peace University, Survey No.\ 124, Paud Road, Kothrud, Pune 411038, Maharashtra, India. E-mail: jyoti.pulgam@mitwpu.edu.in}

\author{Prashant Malavadkar$^2$}
\thanks{$^2$ Corresponding author. Department of Mathematics and Statistics, Dr.\ Vishwanath Karad MIT World Peace University, Survey No.\ 124, Paud Road, Kothrud, Pune 411038, Maharashtra, India. E-mail: prashant.malavadkar@mitwpu.edu.in}


\begin{abstract}
The generalized corona of a graph with a family of graphs, one for each vertex, is obtained by joining every vertex of the base graph to all vertices of its own copy of the corresponding graph; the corona of Frucht and Harary is the case in which all graphs in the family are equal. Each attached copy meets the rest of the graph in a single vertex, so additivity of the orientable genus over blocks reduces the genus of a generalized corona to the genus of the base graph plus the genera of the graphs obtained from the attached graphs by adding a dominating vertex. We study this last quantity through faces that meet every vertex, and we determine it for complete graphs, balanced complete bipartite graphs, wheels, prisms, antiprisms, the Petersen graph (value 2). The Petersen case follows from a face-length criterion. We also obtain bounds under edge addition, a formula for iterated coronas, a criterion for two coronas taken in opposite order to have the same genus, and we show that the corona is not associative with respect to the genus.
\end{abstract}

\maketitle

\section{Introduction}\label{sec:intro}

The orientable genus $\gamma(G)$ of a graph $G$ is the least $g$ such that $G$ can be drawn without crossings on the orientable surface $S_g$ of genus $g$. Exact values are known for few families. Ringel and Youngs \cite{RY68} and Ringel \cite{Ringel65} determined the genus of complete graphs and complete bipartite graphs. The genus of joins and of complete tripartite graphs has been studied in many papers \cite{Craft98,ES07,ES09,ESZ06,Kurauskas17,LvChen19,MaRen19,StahlWhite76,White69}. Computing the genus is NP-hard \cite{Thomassen89}, although for each fixed $g$ there is a linear-time algorithm that decides whether $\gamma(G)\le g$ \cite{Mohar99,KMR08}, and exact programs exist for graphs of moderate size \cite{BCHK16,Brinkmann22}. It is therefore useful to know graph operations for which the genus of the result follows from the genera of the operands. The standard example is additivity over blocks, due to Battle, Harary, Kodama and Youngs \cite{BHKY62}. For amalgamations along two vertices the genus is no longer additive, and only bounds are known \cite{Alpert73,DGH81}.

The corona $G\circ H$ of Frucht and Harary \cite{FruchtHarary70} is obtained from a copy of $G$ and $|V(G)|$ copies of $H$ by joining the $i$-th vertex of $G$ to every vertex of the $i$-th copy of $H$. Iterated coronas have been used as models of large networks \cite{Sharma17}. To our knowledge the orientable genus of coronas has not been studied before. The planarity of coronas is classical, and the structure of a corona as a family of pieces glued at single vertices of $G$ is routine; what is new here is the evaluation of the resulting genus, which depends on the genus of graphs with a dominating vertex.

We work with the \emph{generalized corona}: given a graph $G$ with vertices $v_1,\dots,v_n$ and graphs $H_1,\dots,H_n$, which may all differ, join $v_i$ to every vertex of a copy of $H_i$. We write $G\circ(H_1,\dots,H_n)$ for the result (Definition~\ref{def:gcorona}); the corona is the case $H_1=\dots=H_n=H$. Each copy of $H_i$ meets the rest of the graph in the single vertex $v_i$, which is a cut vertex when $n\ge2$. The blocks of the generalized corona are the blocks of $G$ together with the blocks of the graphs $K_1+H_i$, where $K_1+H$ denotes $H$ with one added vertex adjacent to all vertices of $H$.

In this paper we determine the genus of $G\circ(H_1,\dots,H_n)$ in section \ref{sec:main}.
Its use is that it moves all the difficulty to the genus of $K_1+H$, a graph with a dominating vertex. The key to this quantity is the following observation, proved in Section~\ref{sec:main} as Proposition~\ref{prop:face}.  Section ~\ref{sec:values}, Section~\ref{sec:petersen} and Section ~\ref{sec:faces} determine $\gamma(K_1+H)$  for complete graphs, balanced complete bipartite graphs, wheels, prisms and antiprisms and Petersen graph. Section \ref{sec:bounds} discuss the bounds of \(\gamma(K_1 + H)\) by comparing it with \(\gamma(H)\). Section~\ref{sec:more} gives an exact criterion for $\gamma(G\circ H)=\gamma(H\circ G)$ and shows that the corona is not associative for the genus following with the examples, including an infinite family of coronas of graphs of different genus that have equal genus. Section~\ref{sec:problems} lists open problems.

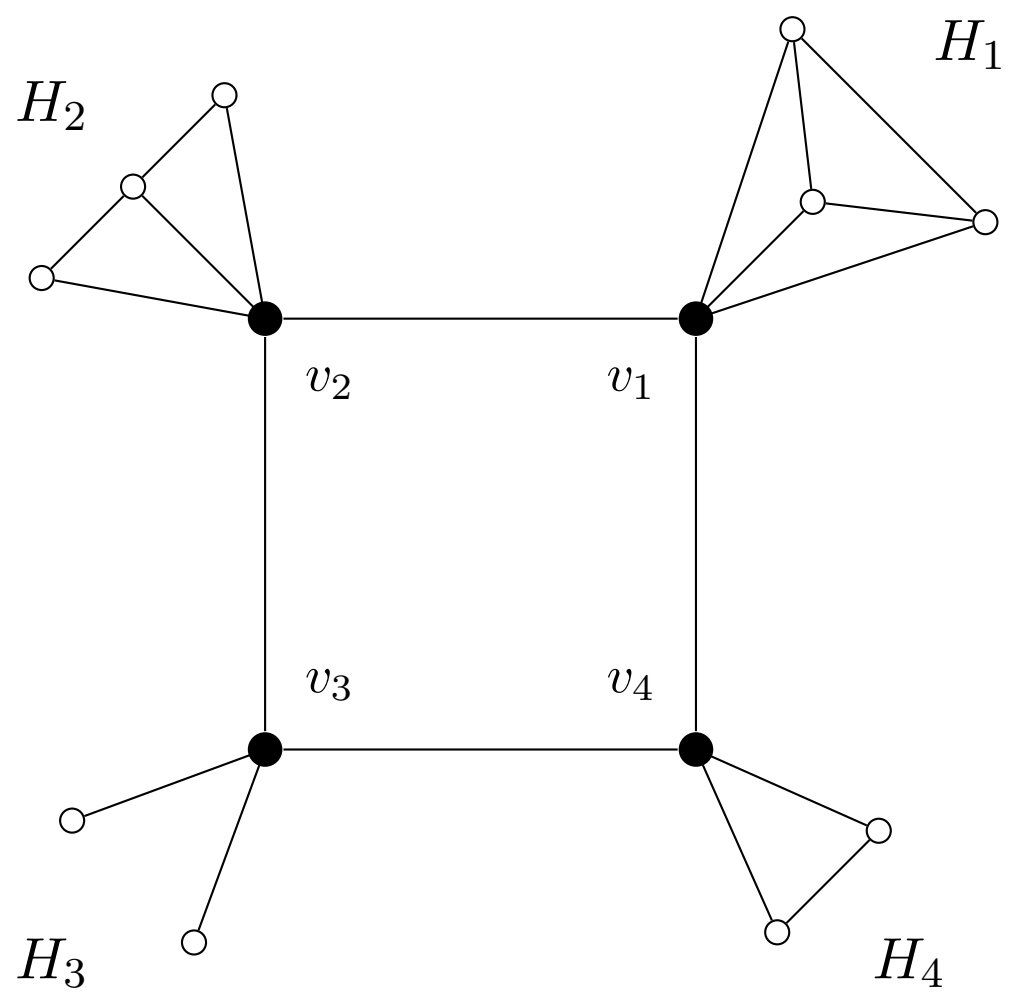
\begin{figure}
\centering
\begin{tikzpicture}[scale=0.95,
  v/.style={circle,fill=black,inner sep=2.3pt},
  w/.style={circle,draw=black,fill=white,inner sep=1.6pt}]
  \begin{scope}[shift={(1.5,1.5)},rotate=45]
    \node[w] (a1) at (1.9,-0.95) {};
    \node[w] (b1) at (1.9,0.95) {};
    \node[w] (c1) at (1.15,0) {};
    \draw (0,0)--(a1); \draw (0,0)--(b1); \draw (0,0)--(c1);
    \draw (a1)--(b1)--(c1)--(a1);
    \node at (2.7,0) {$H_1$};
  \end{scope}
  \begin{scope}[shift={(-1.5,1.5)},rotate=135]
    \node[w] (a2) at (1.3,-0.9) {};
    \node[w] (b2) at (1.3,0) {};
    \node[w] (c2) at (1.3,0.9) {};
    \draw (0,0)--(a2); \draw (0,0)--(b2); \draw (0,0)--(c2);
    \draw (a2)--(b2)--(c2);
    \node at (2.1,0) {$H_2$};
  \end{scope}
  \begin{scope}[shift={(-1.5,-1.5)},rotate=225]
    \node[w] (a3) at (1.3,-0.6) {};
    \node[w] (b3) at (1.3,0.6) {};
    \draw (0,0)--(a3); \draw (0,0)--(b3);
    \node at (2.1,0) {$H_3$};
  \end{scope}
  \begin{scope}[shift={(1.5,-1.5)},rotate=315]
    \node[w] (a4) at (1.3,-0.5) {};
    \node[w] (b4) at (1.3,0.5) {};
    \draw (0,0)--(a4); \draw (0,0)--(b4);
    \draw (a4)--(b4);
    \node at (2.1,0) {$H_4$};
  \end{scope}
  \node[v] (v1) at (1.5,1.5) {};
  \node[v] (v2) at (-1.5,1.5) {};
  \node[v] (v3) at (-1.5,-1.5) {};
  \node[v] (v4) at (1.5,-1.5) {};
  \draw (v1)--(v2)--(v3)--(v4)--(v1);
  \node[font=\small] at (1.05,1.05) {$v_1$};
  \node[font=\small] at (-1.05,1.05) {$v_2$};
  \node[font=\small] at (-1.05,-1.05) {$v_3$};
  \node[font=\small] at (1.05,-1.05) {$v_4$};
\end{tikzpicture}
\caption{The generalized corona $C_4\circ(K_3,P_3,\edgeless{2},K_2)$. The black vertices $v_1,\dots,v_4$ form the cycle $C_4$, and $v_i$ is joined to every vertex of $H_i$, where $H_1=K_3$, $H_2=P_3$, $H_3=\edgeless{2}$ (two isolated vertices) and $H_4=K_2$. If all $H_i$ equal one graph $H$, the result is the corona $C_4\circ H$.}\label{fig:corona}
\end{figure}

\section{Preliminaries}\label{sec:prelim}

All graphs are finite and simple. For a graph $G$ we write $V(G)$ and $E(G)$ for its vertex and edge sets respectively, $N_G(x)$ for the set of neighbours of $x$, and $\beta(G)=|E(G)|-|V(G)|+c(G)$ for the cycle rank, where $c(G)$ is the number of components. The genus $\gamma(G)$ is used for connected graphs. The join $G+H$ of two graphs with disjoint vertex sets is obtained from their disjoint union by adding all edges between $V(G)$ and $V(H)$. So $K_1+H$ is $H$ with a new vertex adjacent to all vertices of $H$, and $K_1+H$ is connected for every graph $H$. A graph is \emph{outerplanar} if it has a plane embedding in which all vertices lie on the boundary of one face.

\subsection{Generalized coronas}

\begin{definition}\label{def:gcorona}\cite{Laali16}
Let $G$ be a graph with $V(G)=\{v_1,\dots,v_n\}$ and let $H_1,\dots,H_n$ be graphs with non-empty vertex sets, pairwise disjoint and disjoint from $V(G)$. The \emph{generalized corona} $G\circ(H_1,\dots,H_n)$ has vertex set $V(G)\cup\bigcup_iV(H_i)$ and edge set
\[
E(G)\ \cup\ \bigcup_{i=1}^nE(H_i)\ \cup\ \bigcup_{i=1}^n\{v_ix:x\in V(H_i)\}.
\]
If $H_i$ is a copy of one graph $H$ for every $i$, this is the \emph{corona} $G\circ H$.
\end{definition}

If $h_i=|V(H_i)|$, then $|V(G\circ(H_1,\dots,H_n))|=n+\sum_ih_i$ and $|E(G\circ(H_1,\dots,H_n))|=|E(G)|+\sum_i(|E(H_i)|+h_i)$. The graph $G$ is a subgraph of $G\circ(H_1,\dots,H_n)$, and the latter is connected if and only if $G$ is connected.

\subsection{Embeddings}
An embedding of a connected graph in a surface is \emph{cellular} if every face is homeomorphic to an open disc. A cellular embedding in an orientable surface is determined up to homeomorphism by a rotation system, that is, a cyclic ordering of the neighbours of each vertex, and every rotation system arises from an embedding \cite{GrossTucker87,MoharThomassen01}. Faces are found by face tracing: from an arc $(a,b)$ the next arc is $(b,c)$, where $c$ follows $a$ in the cyclic order at $b$. The orbits of this rule are the boundary walks of the faces, and the length of a face is the number of arcs in its orbit. Each edge gives two arcs, so the face lengths add up to $2|E|$. A \emph{corner} of a face at a vertex $b$ is a pair of consecutive arcs $(a,b),(b,c)$ of its boundary walk. If $F$ is a face of a cellular embedding, there is a continuous map from the closed unit disc onto the closure of $F$ that is a homeomorphism on the interior and maps the boundary circle onto the boundary walk of $F$, traversing it once \cite{MoharThomassen01}; points of the circle over a vertex $x$ correspond to corners of $F$ at $x$.

\subsection{Standard facts}

\begin{lemma}[Monotonicity]\label{lem:sub}
If $G'$ is a connected subgraph of a connected graph $G$, then $\gamma(G')\le\gamma(G)$.
\end{lemma}

\begin{proof}
An embedding of $G$ in $S_g$ restricts to an embedding of $G'$ in $S_g$.
\end{proof}

\begin{lemma}[Euler's formula; Youngs \cite{Youngs63}]\label{lem:euler}
If a connected graph with $n$ vertices and $e$ edges has a cellular embedding in $S_g$ with $f$ faces, then $n-e+f=2-2g$. If $g=\gamma(G)$, every embedding of $G$ in $S_g$ is cellular.
\end{lemma}

\begin{theorem}[Battle, Harary, Kodama, Youngs \cite{BHKY62}]\label{thm:BHKY}
Let $L$ and $J$ be connected graphs with exactly one common vertex. Then $\gamma(L\cup J)=\gamma(L)+\gamma(J)$.
\end{theorem}

\begin{lemma}[Adding an edge]\label{lem:edge}
If $G$ is a connected graph and $e\notin E(G)$ joins two vertices of $G$, then $\gamma(G)\le\gamma(G+e)\le\gamma(G)+1$.
\end{lemma}

\begin{proof}
The first inequality is Lemma~\ref{lem:sub}. Let $e=uv$ and embed $G$ in $S_g$ with $g=\gamma(G)$; the embedding is cellular by Lemma~\ref{lem:euler}. Choose a face $F_u$ with $u$ on its boundary and a face $F_v$ with $v$ on its boundary. If $F_u=F_v$, draw $e$ inside this face. Otherwise attach a handle to the surface with one end in $F_u$ and the other in $F_v$, and draw $e$ along the handle. The new surface is orientable of genus $g+1$.
\end{proof}

\begin{lemma}[Chartrand and Harary \cite{ChartrandHarary67}]\label{lem:outer}
A graph $H$ is outerplanar if and only if $K_1+H$ is planar. In particular $\gamma(K_1+H)=0$ if and only if $H$ is outerplanar.
\end{lemma}

\begin{proof}
Let $\apex$ be the added vertex. If $K_1+H$ is planar, delete $\apex$ from a plane embedding. All vertices of $H$ lie on the boundary of the face that contained $\apex$, so $H$ is outerplanar. Conversely, if $H$ is embedded in the plane with all vertices on the boundary of one face, place $\apex$ in that face and join it to every vertex of $H$ inside the face.
\end{proof}

\begin{theorem}[Ringel--Youngs \cite{RY68}, Ringel \cite{Ringel65}]\label{thm:RY}
For $p\ge3$, $\gamma(K_p)=\ceil{(p-3)(p-4)/12}$. For $p,q\ge2$, $\gamma(K_{p,q})=\ceil{(p-2)(q-2)/4}$.
\end{theorem}

\begin{theorem}[Kurauskas \cite{Kurauskas17}, Lv and Chen \cite{LvChen19}]\label{thm:tripartite}
For every $p\ge1$, $\gamma(K_{p,p,1})=\ceil{(p-1)(p-2)/4}$. Kurauskas proved the case of even $p$, and Lv and Chen the case of odd $p$.
\end{theorem}

The \emph{Euler genus} $\eg(G)$\cite{StahlBeineke77} of a connected graph is the least $k$ such that $G$ embeds in a surface of Euler characteristic $2-k$, orientable or not. Thus $\eg(G)=\min\{2\gamma(G),\tilde\gamma(G)\}$, where $\tilde\gamma$ is the nonorientable genus, and $\eg(G)\le2\gamma(G)$. We use the nonorientable genus of complete graphs: $\tilde\gamma(K_n)=\ceil{(n-3)(n-4)/6}$ for $n\ge3$, $n\ne7$, and $\tilde\gamma(K_7)=3$ (Ringel \cite{Ringel74}).

\begin{theorem}[Stahl--Beineke \cite{StahlBeineke77}, Miller {\cite[Theorem~1]{Miller87}}]\label{thm:euleradd}
The Euler genus is additive over blocks: if $L$ and $J$ are connected graphs with exactly one common vertex, then $\eg(L\cup J)=\eg(L)+\eg(J)$. It is also additive over edge-amalgamations: if $L$ and $J$ are connected graphs whose intersection is a single edge together with its end vertices, then $\eg(L\cup J)=\eg(L)+\eg(J)$.
\end{theorem}

\section{The genus formula and the face criterion}\label{sec:main}

The formula follows by adding the graphs $K_1+H_i$ to $G$ one at a time; each step is a one-vertex amalgamation, to which Theorem~\ref{thm:BHKY} applies.

\begin{theorem}\label{thm:main}
Let $G$ be a connected graph with $V(G)=\{v_1,\dots,v_n\}$ and let $H_1,\dots,H_n$ be as in Definition~\ref{def:gcorona}. Then
\[
\gamma\bigl(G\circ(H_1,\dots,H_n)\bigr)=\gamma(G)+\sum_{i=1}^n\gamma(K_1+H_i).
\]
\end{theorem}

\begin{proof}
For $i=1,\dots,n$ let $J_i$ be the subgraph of $G\circ(H_1,\dots,H_n)$ with vertex set $\{v_i\}\cup V(H_i)$ and edge set $E(H_i)\cup\{v_ix:x\in V(H_i)\}$. Then $J_i\cong K_1+H_i$, and $J_i$ is connected because $v_i$ is adjacent to every other vertex of $J_i$. Put $L_0=G$ and, for $1\le i\le n$, let $L_i$ be the graph with vertex set $V(L_{i-1})\cup V(J_i)$ and edge set $E(L_{i-1})\cup E(J_i)$.

We show by induction on $i$ that $L_i$ is connected and $\gamma(L_i)=\gamma(G)+\sum_{j\le i}\gamma(K_1+H_j)$. For $i=0$ this holds because $G$ is connected. Let $i\ge1$. By construction $V(L_{i-1})=V(G)\cup V(H_1)\cup\dots\cup V(H_{i-1})$. The sets $V(H_j)$ are pairwise disjoint and disjoint from $V(G)$, and $v_i\in V(G)$, so $V(L_{i-1})\cap V(J_i)=\{v_i\}$. The graphs $L_{i-1}$ (by induction) and $J_i$ are connected, so $L_i$ is connected, and it is the union of two connected graphs with exactly one common vertex. Theorem~\ref{thm:BHKY} gives $\gamma(L_i)=\gamma(L_{i-1})+\gamma(J_i)=\gamma(L_{i-1})+\gamma(K_1+H_i)$, which completes the induction.

Finally $V(L_n)=V(G)\cup\bigcup_iV(H_i)$ and $E(L_n)=E(G)\cup\bigcup_i\bigl(E(H_i)\cup\{v_ix:x\in V(H_i)\}\bigr)$, so $L_n=G\circ(H_1,\dots,H_n)$.
\end{proof}

The graphs $H_i$ need not be connected. The same induction, using Theorem~\ref{thm:euleradd} for blocks, gives $\eg(G\circ(H_1,\dots,H_n))=\eg(G)+\sum_i\eg(K_1+H_i)$.

\begin{corollary}\label{cor:corona}
Let $G$ be a connected graph with $n$ vertices and let $H$ be a graph. Then $\gamma(G\circ H)=\gamma(G)+n\,\gamma(K_1+H)$.
\end{corollary}

\begin{proof}
Take $H_1=\dots=H_n=H$ in Theorem~\ref{thm:main}.
\end{proof}

The next result says exactly when attaching graphs to $G$ does not change the genus or the planarity. Part (b) is a classical planarity criterion for coronas; we include the short proof.

\begin{corollary}\label{cor:outer}
Let $G$ be a connected graph and let $G\circ(H_1,\dots,H_n)$ be a generalized corona. Then
\begin{enumerate}
\item[(a)] $\gamma(G\circ(H_1,\dots,H_n))\ge\gamma(G)$, with equality if and only if every $H_i$ is outerplanar;
\item[(b)] $G\circ(H_1,\dots,H_n)$ is planar if and only if $G$ is planar and every $H_i$ is outerplanar.
\end{enumerate}
\end{corollary}

\begin{proof}
(a) By Theorem~\ref{thm:main}, $\gamma(G\circ(H_1,\dots,H_n))-\gamma(G)=\sum_i\gamma(K_1+H_i)$ is a sum of non-negative integers. It is $0$ if and only if $\gamma(K_1+H_i)=0$ for every $i$, and by Lemma~\ref{lem:outer} this holds if and only if every $H_i$ is outerplanar. For (b), the generalized corona is planar if and only if its genus is $0$, that is, if and only if $\gamma(G)=0$ and every $H_i$ is outerplanar.
\end{proof}

\begin{example}\label{ex:outer}
Paths and cycles are outerplanar, so $\gamma(G\circ P_m)=\gamma(G)$ for $m\ge1$ and $\gamma(G\circ C_m)=\gamma(G)$ for $m\ge3$; for example $\gamma(K_5\circ P_3)=\gamma(K_5)=1$. The graph in Figure~\ref{fig:corona} is planar, since $C_4$ is planar and $K_3$, $P_3$, $\edgeless{2}$ and $K_2$ are outerplanar.
\end{example}

\subsection{The face criterion}
Theorem~\ref{thm:main} reduces everything to $\gamma(K_1+H)$. The next proposition describes this quantity in terms of embeddings of $H$. By a \emph{face} of an embedding of $H$ we mean a connected component of the complement of $H$ in the surface. A face $F$ \emph{sees} a vertex $x$ if $x$ lies in the closure of $F$, and faces $F_1,\dots,F_t$ \emph{cover} $H$ if every vertex of $H$ is seen by at least one of them. A face of a cellular embedding is an open disc, and it sees the vertices on its boundary walk.

\begin{proposition}[Face criterion]\label{prop:face}
Let $H$ be a connected graph and let $g\ge0$.
\begin{enumerate}
\item[(a)] If $K_1+H$ embeds in $S_g$, then $H$ embeds in $S_g$ with a face that sees every vertex of $H$.
\item[(b)] If $H$ has a cellular embedding in $S_g$ and faces $F_1,\dots,F_t$ of this embedding cover $H$, then $K_1+H$ embeds in $S_{g+t-1}$. In particular $\gamma(K_1+H)\le g+t-1$.
\item[(c)] $\gamma(K_1+H)=\gamma(H)$ if and only if $H$ has an embedding in $S_{\gamma(H)}$ with a face that sees every vertex of $H$.
\end{enumerate}
\end{proposition}

\begin{proof}
(a) Delete the added vertex $\apex$ from the embedding. The face of the resulting embedding of $H$ that contains $\apex$ contains the edges $\apex x$ except for their ends $x$, so its closure contains every vertex $x$.

(b) We may assume that $F_1,\dots,F_t$ are distinct and that each $F_j$ with $j\ge2$ sees a vertex not seen by $F_1,\dots,F_{j-1}$ (otherwise discard $F_j$, which only decreases $t$). Let $X_1$ be the set of vertices seen by $F_1$, and for $j\ge2$ let $X_j$ be the set of vertices seen by $F_j$ and by none of $F_1,\dots,F_{j-1}$. For every vertex in $X_j$ fix one corner of $F_j$ at that vertex. Place $\apex$ in the open disc $F_1$ and join $\apex$ to the chosen corners of the vertices of $X_1$ by pairwise disjoint arcs in $F_1$; this is possible because $F_1$ is a disc (use the characteristic map of $F_1$ from Section~\ref{sec:prelim} and draw radial arcs). The arcs divide $F_1$ into sectors. For $j=2,\dots,t$ choose a small closed disc $\Delta_j\subset F_1$ in one of the sectors, close to $\apex$, with the discs $\Delta_j$ pairwise disjoint, and choose a small closed disc $\Delta_j'$ inside $F_j$ disjoint from $H$. Remove the interiors of all $\Delta_j$ and $\Delta_j'$ and, for each $j$, glue a cylinder to the two boundary circles so that the surface stays orientable. Each cylinder joins two points of one connected surface, so it adds one handle, and the new surface is $S_{g+t-1}$. In the sector of $F_1$ containing $\Delta_j$ draw $|X_j|$ disjoint arcs from $\apex$ to $\partial\Delta_j$; continue them along parallel arcs of the $j$-th cylinder to $\partial\Delta_j'$; and finish them in the annulus $F_j\setminus\Delta_j'$ by disjoint radial arcs to the chosen corners of the vertices of $X_j$. Which arc goes to which vertex is ours to choose, so the cyclic order of the endpoints on $\partial\Delta_j'$ can be matched with the order of the corners along the boundary walk of $F_j$. The result is an embedding of $K_1+H$ in $S_{g+t-1}$.

(c) If $\gamma(K_1+H)=\gamma(H)=k$, then (a) gives an embedding of $H$ in $S_k$ with a face that sees every vertex. Conversely, an embedding of $H$ in $S_k$ with $k=\gamma(H)$ is cellular by Lemma~\ref{lem:euler}, so (b) with $t=1$ gives $\gamma(K_1+H)\le k$, and $\gamma(K_1+H)\ge\gamma(H)$ by Lemma~\ref{lem:sub}.
\end{proof}

The case $g=0$, $t=1$ is Lemma~\ref{lem:outer}.

\subsection*{Iterated coronas}
Coronas can be iterated. For graphs $G,H_1,\dots,H_r$ put $G^{(0)}=G$ and $G^{(j)}=G^{(j-1)}\circ H_j$ for $j\ge1$, and let $h_j=|V(H_j)|$ and $n=|V(G)|$. Since $|V(G^{(j)})|=|V(G^{(j-1)})|(1+h_j)$, we have $|V(G^{(j)})|=n\prod_{i=1}^{j}(1+h_i)$.

\begin{theorem}\label{thm:iterated}
Let $G$ be connected and let $H_1,\dots,H_r$ be graphs. Then
\[
\gamma\bigl(G^{(r)}\bigr)=\gamma(G)+n\sum_{i=1}^{r}\Bigl(\prod_{j=1}^{i-1}(1+h_j)\Bigr)\gamma(K_1+H_i),
\]
where an empty product equals $1$.
\end{theorem}

\begin{proof}
We use induction on $r$; the case $r=0$ is trivial. Each $G^{(j)}$ is connected because $G$ is. By Corollary~\ref{cor:corona} applied to $G^{(r-1)}$ and $H_r$,
\begin{align*}
\gamma(G^{(r)})&=\gamma(G^{(r-1)})+|V(G^{(r-1)})|\,\gamma(K_1+H_r)\\
&=\gamma(G^{(r-1)})+n\Bigl(\prod_{j=1}^{r-1}(1+h_j)\Bigr)\gamma(K_1+H_r).
\end{align*}
Substituting the induction hypothesis for $\gamma(G^{(r-1)})$ gives the claim.
\end{proof}

\begin{corollary}\label{cor:powers}
Let $G$ be connected with $n$ vertices, let $H$ be a graph with $h\ge1$ vertices, and let $G\circ^{r}H=(\cdots((G\circ H)\circ H)\cdots)\circ H$ with $r$ coronas. Then
\[
\gamma(G\circ^{r}H)=\gamma(G)+n\,\gamma(K_1+H)\,\frac{(1+h)^r-1}{h}.
\]
\end{corollary}

\begin{proof}
Take $H_1=\dots=H_r=H$ in Theorem~\ref{thm:iterated} and use $\sum_{i=0}^{r-1}(1+h)^i=((1+h)^r-1)/h$.
\end{proof}

\begin{example}\label{ex:iter}
Let $G=K_5$ and $H=K_4$. Then $K_1+K_4=K_5$, so $\gamma(K_1+H)=1$, and $\gamma(K_5\circ^2K_4)=1+5\cdot1\cdot\frac{5^2-1}{4}=31$. Directly, $\gamma(K_5\circ K_4)=1+5=6$; this graph has $25$ vertices, so $\gamma((K_5\circ K_4)\circ K_4)=6+25=31$. The graph $K_5\circ^2K_4$ has $125$ vertices.
\end{example}

\section{Complete and complete bipartite second factors}\label{sec:values}

By Corollary~\ref{cor:corona}, the genus of $G\circ H$ is known as soon as $\gamma(G)$ and $\gamma(K_1+H)$ are known. This section lists cases in which $\gamma(K_1+H)$ follows from known genera.

\begin{proposition}\label{prop:complete}
Let $G$ be a connected graph with $n$ vertices and let $p\ge2$. Then
\[
\gamma(G\circ K_p)=\gamma(G)+n\Bigl\lceil\frac{(p-2)(p-3)}{12}\Bigr\rceil.
\]
In particular $G\circ K_p$ is planar if and only if $G$ is planar and $p\le3$.
\end{proposition}

\begin{proof}
The graph $K_1+K_p$ is $K_{p+1}$ and $p+1\ge3$. By Theorem~\ref{thm:RY}, $\gamma(K_{p+1})=\ceil{(p-2)(p-3)/12}$, which is $0$ if and only if $p\le3$. By using  Corollary~\ref{cor:corona},
\(
\gamma(G\circ K_p)=\gamma(G)+n\Bigl\lceil\frac{(p-2)(p-3)}{12}\Bigr\rceil.
\)
\end{proof}

\begin{proposition}\label{prop:bip}
Let $G$ be a connected graph with $n$ vertices and let $p\ge1$. Then
\[
\gamma(G\circ K_{p,p})=\gamma(G)+n\Bigl\lceil\frac{(p-1)(p-2)}{4}\Bigr\rceil.
\]
\end{proposition}

\begin{proof}
The graph $K_1+K_{p,p}$ is the complete tripartite graph $K_{p,p,1}$. Apply Theorem~\ref{thm:tripartite} and Corollary~\ref{cor:corona}.
\end{proof}
\begin{example}\label{ex:k33}
	Since $K_{3,3}\subset K_{3,3,1}\subset K_7$ and $\gamma(K_{3,3})=\gamma(K_7)=1$, we have $\gamma(K_1+K_{3,3})=1$. Hence $\gamma(K_8\circ K_{3,3})=\gamma(K_8)+8=2+8=10$. Also $\gamma(C_8\circ K_5)=0+8\gamma(K_6)=8$.
\end{example}

For $K_{p,q}$ with $p\ne q$ the genus of $K_{p,q,1}$ is not known in general. White's conjecture \cite{White69} states that $\gamma(K_{l,m,n})=\ceil{(l-2)(m+n-2)/4}$ whenever $l\ge m\ge n$; the right side is the genus of $K_{l,m+n}$, which is a subgraph of $K_{l,m,n}$. The conjecture is open in general. It has been confirmed when all part sizes are even and for $K_{l,m,m}$ with $l\ge m\ge2$ (see the introduction of \cite{Kurauskas17}), and for $K_{p,p,1}$ by Theorem~\ref{thm:tripartite}. The analogous conjecture of Stahl and White \cite{StahlWhite76} for the nonorientable genus has been settled, with three exceptions \cite{ESZ06}.

\begin{proposition}\label{prop:unbal}
Let $p\ge q\ge1$ with $p\ge2$, and let $G$ be connected with $n$ vertices. Then
\[
\gamma(G)+n\Bigl\lceil\frac{(p-2)(q-1)}{4}\Bigr\rceil\ \le\ \gamma(G\circ K_{p,q})\ \le\ \gamma(G)+n\Bigl\lceil\frac{(p-1)(p-2)}{4}\Bigr\rceil.
\]
\end{proposition}

\begin{proof}
The graph $K_1+K_{p,q}=K_{p,q,1}$ contains $K_{p,q+1}$ as a subgraph, so by Lemma~\ref{lem:sub} and Theorem~\ref{thm:RY}, $\gamma(K_{p,q,1})\ge\ceil{(p-2)(q-1)/4}$. It is a subgraph of $K_{p,p,1}$, so $\gamma(K_{p,q,1})\le\ceil{(p-1)(p-2)/4}$ by Theorem~\ref{thm:tripartite}. Apply Corollary~\ref{cor:corona}.
\end{proof}

\begin{corollary}\label{cor:unbal1}
For $p\ge2$ we have $K_1\circ K_{p,p+1}=K_{p,p+1,1}$ and
\[
\left\lceil\frac{(p-1)^2}{4}\right\rceil\ \le\ \gamma(K_{p,p+1,1})\ \le\ \left\lceil\frac{p(p-1)}{4}\right\rceil .
\]
The two bounds coincide if and only if $p\in\{2,4\}$. Consequently $\gamma(K_{2,3,1})=1$ and $\gamma(K_{4,5,1})=3$.
\end{corollary}

\begin{proof}
Apply Proposition~\ref{prop:unbal} with $G=K_1$, $n=1$, $\gamma(G)=0$, and with $(p+1,p)$ in place of $(p,q)$. We evaluate the two bounds according to the residue of $p$ modulo $4$. For $p=4k$ they are $4k^2-2k+1$ and $4k^2-k$; for $p=4k+1$ they are $4k^2$ and $4k^2+k$; for $p=4k+2$ they are $4k^2+2k+1$ and $4k^2+3k+1$; for $p=4k+3$ they are $4k^2+4k+1$ and $4k^2+5k+2$. The differences, upper minus lower, are $k-1$, $k$, $k$ and $k+1$ respectively. For $p\ge2$ they vanish only for $p=4$ (in the first case with $k=1$) and $p=2$ (in the third case with $k=0$). The common values are $1$ for $p=2$ and $3$ for $p=4$.
\end{proof}

White's conjecture predicts $\gamma(K_{p+1,p,1})=\lceil(p-1)^2/4\rceil$, that is, the lower bound of Corollary~\ref{cor:unbal1}.

\begin{example}\label{ex:k45}
Let $G$ be a connected graph with $n$ vertices. By Corollaries~\ref{cor:corona} and~\ref{cor:unbal1},
\begin{enumerate}
\item[\textnormal{(i)}] $\gamma(G\circ K_{2,3})=\gamma(G)+n$,
\item[\textnormal{(ii)}] $\gamma(G\circ K_{4,5})=\gamma(G)+3n$.
\end{enumerate}
\end{example}

\begin{corollary}\label{cor:unbal2}
Let $G$ be a connected graph on $n$ vertices and let $p\ge1$. The genus of $G\circ K_{p,p+1}$ is bounded as follows, according to the residue class of $p$ modulo $4$ (here $k\ge0$, and $k\ge1$ in (i)):
\begin{enumerate}
\item[\textnormal{(i)}] if $p=4k$, then $\gamma(G)+n(4k^2-2k+1)\le\gamma(G\circ K_{4k,4k+1})\le\gamma(G)+n(4k^2-k)$;
\item[\textnormal{(ii)}] if $p=4k+1$, then $\gamma(G)+4nk^2\le\gamma(G\circ K_{4k+1,4k+2})\le\gamma(G)+n(4k^2+k)$;
\item[\textnormal{(iii)}] if $p=4k+2$, then $\gamma(G)+n(4k^2+2k+1)\le\gamma(G\circ K_{4k+2,4k+3})\le\gamma(G)+n(4k^2+3k+1)$;
\item[\textnormal{(iv)}] if $p=4k+3$, then $\gamma(G)+n(4k^2+4k+1)\le\gamma(G\circ K_{4k+3,4k+4})\le\gamma(G)+n(4k^2+5k+2)$.
\end{enumerate}
\end{corollary}

\begin{proof}
Apply Corollary~\ref{cor:corona} and the bounds of Corollary~\ref{cor:unbal1}, evaluated as in its proof.
\end{proof}

\subsection*{Familiar graphs as coronas}
For $G=K_1$ the corona $K_1\circ H$ is the join $K_1+H$, so several familiar graphs are coronas. Table~\ref{tab:standard} lists the cases that follow from the results above. For the windmill graph, $K_1+(rK_k)$ is the union of $r$ copies of $K_{k+1}$ sharing one vertex, and Theorem~\ref{thm:BHKY} gives genus $r\,\gamma(K_{k+1})$. For a connected graph $G$, the corona $G\circ K_1$ is $G$ with a pendant vertex attached to each vertex (for $G=C_n$ it is the sun graph), and $G\circ\edgeless{r}$ has $r$ pendant vertices at each vertex; both have the genus of $G$ by Corollary~\ref{cor:outer}.

\begin{table}[ht]
\centering
\caption{Familiar graphs written as coronas.}\label{tab:standard}
\begin{tabular}{@{}llc@{}}
\toprule
graph & corona form & genus\\
\midrule
wheel with $m\ge3$ rim vertices & $K_1\circ C_m$ & $0$\\
fan on $m$ path vertices & $K_1\circ P_m$ & $0$\\
complete graph $K_{p+1}$ & $K_1\circ K_p$ & $\ceil{(p-2)(p-3)/12}$\\
complete tripartite graph $K_{p,p,1}$ & $K_1\circ K_{p,p}$ & $\ceil{(p-1)(p-2)/4}$\\
sun graph on $2n$ vertices & $C_n\circ K_1$ & $0$\\
windmill ($r$ copies of $K_{k+1}$ at a vertex) & $K_1\circ(rK_k)$ & $r\ceil{(k-2)(k-3)/12}$\\
\bottomrule
\end{tabular}
\end{table}

\subsection*{The increment for complete and complete bipartite graphs}
For a connected graph $H$ put $\delta(H)=\gamma(K_1+H)-\gamma(H)$. By Lemma~\ref{lem:sub}, $\delta(H)\ge0$. Then $\gamma(G\circ H)=\gamma(G)+|V(G)|(\gamma(H)+\delta(H))$, so $\delta(H)$ is the extra genus per vertex of $G$ caused by the join with the hub. Theorems~\ref{thm:RY} and~\ref{thm:tripartite} give $\delta$ exactly for two infinite families of non-planar graphs.

\begin{proposition}\label{prop:deltaK}
For $p\ge3$, $\delta(K_p)=\ceil{(p-2)(p-3)/12}-\ceil{(p-3)(p-4)/12}$, and for $p\ge2$, $\delta(K_{p,p})=\ceil{(p-1)(p-2)/4}-\ceil{(p-2)^2/4}$. Both quantities are unbounded as $p\to\infty$ (unboundedness of $\delta$ already occurs on planar graphs, see the remark after Proposition~\ref{prop:eulerbound}).
\end{proposition}

\begin{proof}
The formulas follow from $K_1+K_p=K_{p+1}$, $K_1+K_{p,p}=K_{p,p,1}$ and Theorems~\ref{thm:RY} and~\ref{thm:tripartite}. Since $\ceil{x}-\ceil{y}>x-y-1$, we get $\delta(K_p)>(p-3)/6-1$ and $\delta(K_{p,p})>(p-2)/4-1$.
\end{proof}

For example $\delta(K_4)=1$, $\delta(K_5)=\delta(K_6)=0$, $\delta(K_7)=1$, $\delta(K_{3,3})=0$ and $\delta(K_{4,4})=1$.

\section{A face-length criterion and the Petersen graph}\label{sec:petersen}

For complete and complete bipartite graphs, $\gamma(K_1+H)$ is known from the literature. For sparse graphs $H$ the Euler bound of Proposition~\ref{prop:eulerbound} below is weak, and a different argument is needed. We give one that uses the girth and the hamiltonicity of $H$, and apply it first to the Petersen graph.

\begin{lemma}\label{lem:facelength}
Let $G$ be a connected graph of minimum degree at least $2$ and girth $g$, with a cellular embedding in an orientable surface. Then every face has length at least $g$.
\end{lemma}

\begin{proof}
Let $W$ be the boundary walk of a face, of length $\ell\ge1$, and let $T$ be the subgraph formed by the edges of $W$. Then $T$ is connected. Suppose $T$ has no cycle, so $T$ is a tree with at least one edge. Fix a vertex $s$ of $T$ and let $u$ be a vertex of $T$ at maximum distance from $s$. Then $u$ is a leaf of $T$, with unique edge $xu$ in $T$. Each time $W$ reaches $u$ it arrives along $xu$, and the next edge of $W$ is $uy$, where $y$ follows $x$ in the cyclic order at $u$. Since $uy\in E(T)$, we get $y=x$. So $x$ is the only neighbour of $u$ in $G$, which contradicts the minimum degree. Hence $T$ contains a cycle, which has at least $g$ edges. Since $W$ traverses each edge of $T$ at least once, $\ell\ge g$.
\end{proof}

\begin{lemma}\label{lem:nonham}
A cubic graph with $10$ vertices and girth at least $5$ is not hamiltonian.
\end{lemma}

\begin{proof}
Suppose such a graph has a hamiltonian cycle $C=c_0c_1\cdots c_9c_0$, with indices modulo $10$. Each vertex has degree $3$, so it lies on exactly one edge not in $C$; call these five edges chords. Let $\{c_i,c_j\}$ be a chord and let $d\in\{2,3,4,5\}$ be the distance between $i$ and $j$ modulo $10$. The chord and the shorter arc of $C$ form a cycle of length $d+1$, so $d+1\ge5$ and $d\in\{4,5\}$.

Suppose some chord has $d=4$. After renaming the vertices of $C$ we may assume it is $\{c_0,c_4\}$. The chord at $c_9$ ends in a vertex at distance at least $4$ from $c_9$, that is, in $c_3$, $c_4$ or $c_5$. It is not $c_4$, which already lies on the chord $\{c_0,c_4\}$. If it ends in $c_3$, then $c_9c_0c_4c_3$ is a $4$-cycle, and if it ends in $c_5$, then $c_9c_0c_4c_5$ is a $4$-cycle. Both contradict the girth. Hence every chord has $d=5$, so the chords are $\{c_i,c_{i+5}\}$ for $i=0,\dots,4$. But then $c_0c_1c_6c_5$ is a $4$-cycle, a contradiction.
\end{proof}

\begin{figure}[ht]
\centering
\begin{tikzpicture}[scale=1.15,
	v/.style={circle,draw=black,fill=white,inner sep=1.5pt,font=\small}]
	
	\foreach \i in {0,...,4}{
		\node[v] (o\i) at ({90+72*\i}:2) {\i};
		\pgfmathtruncatemacro{\j}{\i+5}
		\node[v] (i\i) at ({90+72*\i}:1) {\j};
	}
	
	\foreach \i/\j in {0/1,1/2,2/3,3/4,4/0}{
		\draw (o\i)--(o\j);
	}
	
	\draw (o0)--(i0); 
	\foreach \i in {1,...,4}{
		\draw (o\i)--(i\i);
	}
	
	\foreach \i/\j in {0/2,2/4,4/1,1/3,3/0}{
		\draw (i\i)--(i\j);
	}
	
\end{tikzpicture}
\caption{The Petersen graph $\Pet$ with the labelling used in the text. The outer cycle is $0,1,2,3,4$, the spokes are $\{i,i+5\}$, and the inner edges form the pentagram $5,7,9,6,8$.}\label{fig:petersen}
\end{figure}
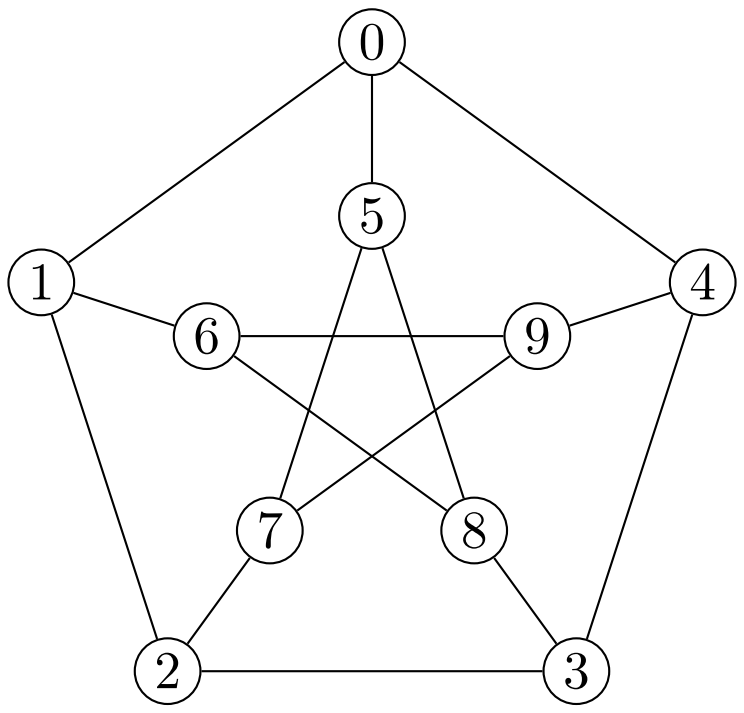

Let $\Pet$ be the Petersen graph with vertices $0,1,\dots,9$ (Figure~\ref{fig:petersen}). Its edges are those of the outer cycle $0,1,2,3,4$, the spokes $\{i,i+5\}$ for $i=0,\dots,4$, and the inner edges $\{i+5,\,5+((i+2)\bmod5)\}$ for $i=0,\dots,4$. It is cubic, with $10$ vertices and $15$ edges.

\begin{lemma}\label{lem:petersen}
The Petersen graph has girth $5$, is not planar, is not hamiltonian, and $\gamma(\Pet)=1$.
\end{lemma}

\begin{proof}
\emph{Girth.} The neighbourhoods are $N(0)=\{1,4,5\}$, $N(1)=\{0,2,6\}$, $N(2)=\{1,3,7\}$, $N(3)=\{2,4,8\}$, $N(4)=\{0,3,9\}$, $N(5)=\{0,7,8\}$, $N(6)=\{1,8,9\}$, $N(7)=\{2,5,9\}$, $N(8)=\{3,5,6\}$, $N(9)=\{4,6,7\}$. Inspection shows that adjacent vertices have no common neighbour and that two distinct vertices have at most one common neighbour. A triangle would give two adjacent vertices with a common neighbour, and a $4$-cycle would give two vertices with two common neighbours. So $\Pet$ has no cycle of length $3$ or $4$, and the cycle $0,1,2,3,4$ has length $5$.

\emph{Non-planarity.} A plane embedding of $\Pet$ would have $f=2-10+15=7$ faces, each of length at least $5$ by Lemma~\ref{lem:facelength}, and then $5\cdot7\le2\cdot15$, which is false.

\emph{Genus one.} The rotation system in Table~\ref{tab:rotP} has, by face tracing, the five faces $(0,4,3,8,6,1)$, $(0,5,7,9,4)$, $(0,1,2,3,4,9,6,8,5)$, $(1,6,9,7,2)$ and $(2,7,5,8,3)$ of lengths $6,5,9,5,5$. These add up to $30=2\cdot15$, and $10-15+5=0$, so this is a cellular embedding in the torus.

\emph{Non-hamiltonicity.} This is Lemma~\ref{lem:nonham}.
\end{proof}

When $K_1+H$ has the same genus as $H$, the new vertex lies in one face of a minimum genus embedding of $H$, and this face meets every vertex of $H$ (Proposition~\ref{prop:face}(a)). If $H$ has large girth, such a face is too long for the total face length available. The following proposition makes this precise.

\begin{proposition}\label{prop:criterion}
Let $H$ be a connected graph with $\nu\ge3$ vertices, $\varepsilon$ edges, minimum degree at least $2$ and girth $g$. Put $k=\gamma(H)$ and $f_k=\varepsilon-\nu+2-2k$. If
\begin{enumerate}
\item[(i)] $\nu+(f_k-1)g>2\varepsilon$, or
\item[(ii)] $\nu+(f_k-1)g=2\varepsilon$ and $H$ is not hamiltonian,
\end{enumerate}
then $\gamma(K_1+H)\ge\gamma(H)+1$.
\end{proposition}

\begin{proof}
Suppose $K_1+H$ embeds in $S_k$. By Proposition~\ref{prop:face}(a), $H$ embeds in $S_k$ with a face $F$ whose closure contains every vertex of $H$. Since $k=\gamma(H)$, this embedding is cellular by Lemma~\ref{lem:euler}, and by Euler's formula it has $f_k$ faces, whose lengths add up to $2\varepsilon$. Every vertex of $H$ lies on the boundary walk of $F$, so $\ell(F)\ge\nu$. By Lemma~\ref{lem:facelength} every other face has length at least $g$. Therefore
\[
2\varepsilon\ \ge\ \ell(F)+(f_k-1)g\ \ge\ \nu+(f_k-1)g.
\]
This contradicts (i). Under (ii) all inequalities are equalities, so $\ell(F)=\nu$, and a closed walk of length $\nu$ through all $\nu$ vertices passes through each vertex exactly once; since $\nu\ge3$ it is a hamiltonian cycle of $H$, which contradicts (ii). So $K_1+H$ does not embed in $S_k$. An embedding in $S_j$ with $j<k$ would give one in $S_k$ by adding handles, so $K_1+H$ embeds in no $S_j$ with $j\le k$, that is, $\gamma(K_1+H)\ge k+1$.
\end{proof}

\begin{theorem}\label{thm:petersen}
$\gamma(K_1+\Pet)=2$.
\end{theorem}

\begin{proof}
\emph{Lower bound.} For $H=\Pet$ we have $\nu=10$, $\varepsilon=15$, $g=5$, $k=1$ and $f_1=15-10+2-2=5$, so $\nu+(f_1-1)g=10+20=30=2\varepsilon$. Since $\Pet$ is not hamiltonian (Lemma~\ref{lem:petersen}), Proposition~\ref{prop:criterion}(ii) gives $\gamma(K_1+\Pet)\ge2$.

\emph{Upper bound.} Write $\apex$ for the added vertex. The rotation system in Table~\ref{tab:rotKP} has, at each vertex, exactly the neighbours of that vertex in $K_1+\Pet$. Face tracing gives the twelve faces in Table~\ref{tab:facesKP}, of lengths $3$ (eight faces), $5$ (two faces), $6$ and $10$. These add up to $50=2\cdot25$, as they must. The graph $K_1+\Pet$ has $11$ vertices and $25$ edges, so by Euler's formula the surface has genus $(2-11+25-12)/2=2$. Hence $\gamma(K_1+\Pet)\le2$.
\end{proof}

\begin{corollary}\label{cor:petersen}
For every connected graph $G$ with $n$ vertices, $\gamma(G\circ \Pet)=\gamma(G)+2n$. For example $\gamma(K_5\circ \Pet)=1+10=11$ and $\gamma(\Pet\circ K_{3,3})=\gamma(\Pet)+10\,\gamma(K_{3,3,1})=1+10=11$.
\end{corollary}

\begin{proof}
The first formula is Corollary~\ref{cor:corona} with Theorem~\ref{thm:petersen}. For $\Pet\circ K_{3,3}$ use $\gamma(\Pet)=1$ and $\gamma(K_{3,3,1})=1$ (Example~\ref{ex:k33}).
\end{proof}

The Petersen graph is the extreme case of the following statement.

\begin{lemma}\label{lem:moore}
A cubic graph of girth at least $5$ has at least $10$ vertices, and a cubic graph of girth at least $6$ has at least $14$ vertices.
\end{lemma}

\begin{proof}
Let the girth be at least $5$ and let $x$ be a vertex. The three neighbours of $x$ and the six further neighbours of these (two for each) are all distinct and distinct from $x$: a coincidence or an edge among them would give a cycle of length at most $4$. So there are at least $1+3+6=10$ vertices. Let the girth be at least $6$ and let $xy$ be an edge. The two other neighbours of $x$, the two other neighbours of $y$, and the eight further neighbours of these four vertices are distinct from each other and from $x,y$: a coincidence would give a cycle of length at most $5$. So there are at least $2+4+8=14$ vertices.
\end{proof}

\begin{corollary}\label{cor:cubic}
Let $H$ be a cubic graph with $\gamma(H)=1$ and girth at least $5$. Then $\gamma(K_1+H)\ge2$, and $\gamma(G\circ H)\ge\gamma(G)+2|V(G)|$ for every connected graph $G$.
\end{corollary}

\begin{proof}
Let $H$ have $\nu$ vertices and girth $g\ge5$. Then $\varepsilon=3\nu/2$, $k=1$ and $f_1=\nu/2$, so $\nu+(f_1-1)g-2\varepsilon=\nu+(\nu/2-1)g-3\nu=\tfrac12\bigl(\nu(g-4)-2g\bigr)$. Condition (i) of Proposition~\ref{prop:criterion} is therefore equivalent to $\nu(g-4)>2g$.

Let $g=5$. By Lemma~\ref{lem:moore}, $\nu\ge10$. If $\nu>10$, then $\nu(g-4)=\nu>10=2g$ and (i) holds. If $\nu=10$, then $\nu+(f_1-1)g=2\varepsilon$, and $H$ is not hamiltonian by Lemma~\ref{lem:nonham}, so (ii) holds.

Let $g\ge6$. By Lemma~\ref{lem:moore}, $\nu\ge14$, so $\nu(g-4)\ge14(g-4)=2g+(12g-56)>2g$ and (i) holds.

In all cases Proposition~\ref{prop:criterion} gives $\gamma(K_1+H)\ge2$. The last statement follows from Corollary~\ref{cor:corona}.
\end{proof}

\section{Planar graphs covered by few faces}\label{sec:faces}

For a plane graph $H$ and a face $F$, let $c(F)$ be the least number of further faces which, together with $F$, cover $H$ in the sense of Section~\ref{sec:main}.

\begin{proposition}\label{prop:cover}
Let $H$ be a connected plane graph. Then
\[
\gamma(K_1+H)\ \le\ \min_F c(F),
\]
where the minimum is over all faces $F$ of the given embedding. In particular, if $F_1,\dots,F_c$ are distinct faces that cover $H$, then $\gamma(K_1+H)\le c-1$. If two faces cover $H$ and $H$ has minimum degree at least $3$, then $\gamma(K_1+H)=1$ and $\delta(H)=1$.
\end{proposition}

\begin{proof}
A connected plane graph is a cellular embedding in $S_0$. If $F$ and $c(F)$ further faces cover $H$, then Proposition~\ref{prop:face}(b) with $g=0$ and $t=c(F)+1$ gives $\gamma(K_1+H)\le c(F)$. If $F_1,\dots,F_c$ cover $H$, then $c(F_1)\le c-1$. Now let two faces cover $H$ and let $H$ have minimum degree at least $3$. Then $\gamma(K_1+H)\le1$. Since an outerplanar graph has a vertex of degree at most $2$, $H$ is not outerplanar, so $\gamma(K_1+H)\ge1$ by Lemma~\ref{lem:outer}. Hence $\gamma(K_1+H)=1$, and $\delta(H)=1-\gamma(H)=1$ because $H$ is planar.
\end{proof}

For $c(F)=0$ the proposition is the easy half of Lemma~\ref{lem:outer}.

\begin{corollary}\label{cor:families}
Let $m\ge3$. For the wheel $W_m=K_1+C_m$, the prism $C_m\square K_2$ and the antiprism on $2m$ vertices, $\gamma(K_1+H)=1$ and $\delta(H)=1$. Hence $\gamma(G\circ H)=\gamma(G)+|V(G)|$ for every connected graph $G$ and each such $H$. Moreover $\gamma(K_2+C_m)=1$.
\end{corollary}

\begin{proof}
Each of these graphs is planar with minimum degree at least $3$. For the prism and the antiprism, the two $m$-gonal faces together contain all vertices. For the wheel, the outer face contains the rim and any triangular face contains the hub; the rim vertices have degree $3$ and the hub has degree $m$. Apply Proposition~\ref{prop:cover}. Finally $K_2+C_m=K_1+W_m$.
\end{proof}

The join $C_m+K_n$ of a cycle with a complete graph is treated in \cite{MaRen19} for $n\ge4$ and $m$ large; the case $n=2$ in Corollary~\ref{cor:families} lies outside that range. The prism for $m=4$ is the cube, and the antiprism for $m=3$ is the octahedron. For the fan $F_m=K_1+P_m$ one has $\delta(F_m)=0$, because $F_m$ is outerplanar.

The bound $c-1$ is not sharp in general. A check of all subsets of the twelve faces of the dodecahedron shows that at least $6$ faces are needed to cover its $20$ vertices, so Proposition~\ref{prop:cover} gives $\gamma(K_1+H)\le5$, but $\gamma(K_1+H)=3$ (Table~\ref{tab:computed}).


\begin{table}[ht]
\centering
\caption{Controls for the exhaustive search: number of labellings (without pruning) that give a rotation system of $K_1+H$ with the target number of faces; all programs give the same count.}\label{tab:controls}
\begin{tabular}{@{}lccc@{}}
\toprule
$H$ & $\gamma(K_1+H)$ & target faces & labellings found\\
\midrule
$K_4$ & $1$ & $5$ & $294$\\
$K_{3,3}$ & $1$ & $8$ & $588$\\
triangular prism & $1$ & $8$ & $662$\\
cube $Q_3$ & $1$ & $11$ & $108$\\
Petersen graph & $2$ & $12$ & $148{,}896$\\
\bottomrule
\end{tabular}
\end{table}

\begin{table}[ht]
\centering
\caption{Values of $\gamma(H)$, $\gamma(K_1+H)$ and $\delta(H)$. Entries marked $^{*}$ are proved in this paper; others were computed with the exact genus program \texttt{multi\_genus} \cite{Brinkmann22} (version~0.1.3 of the Python package \texttt{graph-genus}) and are not proved here; for the icosahedron the bound $\gamma(K_1+H)\ge2$ follows from Proposition~\ref{prop:eulerbound} and for the M\"obius--Kantor graph from Corollary~\ref{cor:cubic}, and in all other unstarred entries both bounds rest on the program. $GP(n,k)$ is the generalized Petersen graph with outer cycle $u_0,\dots,u_{n-1}$, inner vertices $w_0,\dots,w_{n-1}$ and edges $u_iu_{i+1}$, $u_iw_i$, $w_iw_{i+k}$ (indices modulo $n$); $GP(5,2)$ is the Petersen graph, $GP(6,2)$ the D\"urer graph, $GP(8,3)$ the M\"obius--Kantor graph, $GP(10,2)$ the dodecahedron, $GP(10,3)$ the Desargues graph and $GP(12,5)$ the Nauru graph. $Q_n$ is the $n$-cube.}\label{tab:computed}
\small
\begin{tabular}{@{}lcccc@{}}
\toprule
$H$ & $|V(H)|$ & $\gamma(H)$ & $\gamma(K_1+H)$ & $\delta(H)$\\
\midrule
cube $Q_3$ & $8$ & $0$ & $1^{*}$ & $1$\\
Petersen $GP(5,2)$ & $10$ & $1$ & $2^{*}$ & $1$\\
icosahedron & $12$ & $0$ & $2$ & $2$\\
$GP(6,2)$ & $12$ & $0$ & $2$ & $2$\\
$GP(7,2)$ & $14$ & $1$ & $2$ & $1$\\

$GP(8,2)$ & $16$ & $0$ & $2$ & $2$\\
M\"obius--Kantor $GP(8,3)$ & $16$ & $1$ & $3$ & $2$\\
hypercube $Q_4$ & $16$ & $1$ & $3$ & $2$\\
$GP(9,2)$ & $18$ & $1$ & $3$ & $2$\\
dodecahedron $GP(10,2)$ & $20$ & $0$ & $3$ & $3$\\
Desargues $GP(10,3)$ & $20$ & $2$ & $3$ & $1$\\
Nauru $GP(12,5)$ & $24$ & $1$ & $4$ & $3$\\
\bottomrule
\end{tabular}
\end{table}

Table~\ref{tab:computed} shows computed values for further graphs. The M\"obius--Kantor graph is another cubic toroidal graph of girth $6$, for which Corollary~\ref{cor:cubic} gives only $\gamma(K_1+H)\ge2$. The dodecahedron shows that $\delta(H)$ can exceed $\gamma(H)$.

\section{Bounds from subgraphs and edge addition}\label{sec:bounds}

When $\gamma(K_1+H)$ is not known, it can be bounded by comparing $H$ with related graphs.

\begin{proposition}[Monotonicity]\label{prop:mono}
Let $G_1\subseteq G_2$ be connected graphs and let $H_1\subseteq H_2$ be graphs with non-empty vertex sets. Then $\gamma(G_1\circ H_1)\le\gamma(G_2\circ H_2)$.
\end{proposition}

\begin{proof}
Label the vertices of $G_2$ so that $V(G_1)\subseteq V(G_2)$. For a vertex $x$ of $G_1$, the copy of $H_1$ attached to $x$ in $G_1\circ H_1$ is a subgraph of the copy of $H_2$ attached to $x$ in $G_2\circ H_2$. Hence $G_1\circ H_1$ is isomorphic to a connected subgraph of $G_2\circ H_2$, and Lemma~\ref{lem:sub} applies.
\end{proof}

\begin{proposition}[Edge addition]\label{prop:edge}
Let $G$ be connected with $n$ vertices, let $H$ be a graph, and let $H'=H+e$ be obtained by adding one edge $e$ between two vertices of $H$. Then
\[
\gamma(G\circ H)\ \le\ \gamma(G\circ H')\ \le\ \gamma(G\circ H)+n.
\]
\end{proposition}

\begin{proof}
The first inequality is Proposition~\ref{prop:mono}. Since $K_1+H'=(K_1+H)+e$, Lemma~\ref{lem:edge} gives $\gamma(K_1+H)\le\gamma(K_1+H')\le\gamma(K_1+H)+1$, and Corollary~\ref{cor:corona} gives the second inequality.
\end{proof}

\begin{corollary}\label{cor:cyclerank}
Let $H$ be connected. Then $\gamma(K_1+H)\le\beta(H)$, and for every connected $G$ with $n$ vertices
\[
\gamma(G)\ \le\ \gamma(G\circ H)\ \le\ \gamma(G)+n\,\beta(H).
\]
In particular $0\le\delta(H)\le\beta(H)$.
\end{corollary}

\begin{proof}
A spanning tree $T$ of $H$ is outerplanar, so $K_1+T$ is planar by Lemma~\ref{lem:outer}. The graph $K_1+H$ is obtained from $K_1+T$ by adding $\beta(H)$ edges, so Lemma~\ref{lem:edge} applied $\beta(H)$ times gives $\gamma(K_1+H)\le\beta(H)$. The upper bound for $\gamma(G\circ H)$ follows from Corollary~\ref{cor:corona}, and the lower bound from Corollary~\ref{cor:outer}. Finally $\delta(H)=\gamma(K_1+H)-\gamma(H)\le\beta(H)$.
\end{proof}

The upper bound is far from sharp for dense graphs: for $H=K_p$ it gives $n(p-1)(p-2)/2$, whereas Proposition~\ref{prop:complete} gives about $np^2/12$. A lower bound comes from Euler's formula.

\begin{proposition}\label{prop:eulerbound}
Let $G$ be connected with $n$ vertices and let $H$ be a graph with $\nu\ge2$ vertices and $\varepsilon$ edges. Then
\[
\gamma(G\circ H)\ \ge\ \gamma(G)+n\Bigl\lceil\frac{\varepsilon-2\nu+3}{6}\Bigr\rceil.
\]
\end{proposition}

\begin{proof}
The graph $K_1+H$ is connected and simple with $\nu+1\ge3$ vertices and $\varepsilon+\nu$ edges. Take a minimum genus embedding, which is cellular by Lemma~\ref{lem:euler}, and let $f$ be its number of faces. Every face has length at least $3$, so $3f\le2(\varepsilon+\nu)$. Euler's formula gives $2-2\gamma(K_1+H)=(\nu+1)-(\varepsilon+\nu)+f=1-\varepsilon+f$, so $f=1+\varepsilon-2\gamma(K_1+H)$. Then $3+3\varepsilon-6\gamma(K_1+H)\le2\varepsilon+2\nu$, that is $\gamma(K_1+H)\ge(\varepsilon-2\nu+3)/6$. Apply Corollary~\ref{cor:corona}.
\end{proof}

For $H=K_p$ this bound equals the Ringel--Youngs value in Proposition~\ref{prop:complete}. If $H$ is a planar triangulation on $\nu$ vertices, then $\varepsilon=3\nu-6$ and the bound gives $\gamma(K_1+H)\ge\lceil(\nu-3)/6\rceil$ while $\gamma(H)=0$, so $\delta$ is unbounded already on planar graphs. For sparse graphs it is weak; for the Petersen graph it gives only $\gamma(K_1+\Pet)\ge0$, which is why Theorem~\ref{thm:petersen} needs the argument of Section~\ref{sec:petersen}.


\section{Further consequences and examples}\label{sec:more}

\subsection{Symmetry and associativity}
The corona is not commutative, and $G\circ H$ and $H\circ G$ usually have different orders. Corollary~\ref{cor:corona} decides exactly when they have the same genus; this is a natural question because both graphs are built from the same two operands. Recall that $\delta(H)=\gamma(K_1+H)-\gamma(H)$ for a connected graph $H$.

\begin{theorem}\label{thm:sym}
Let $G$ and $H$ be connected graphs with $|V(G)|=n$ and $|V(H)|=\mu$. Then $\gamma(G\circ H)=\gamma(H\circ G)$ if and only if
\[
(n-1)\gamma(H)-(\mu-1)\gamma(G)=\mu\,\delta(G)-n\,\delta(H).
\]
\end{theorem}

\begin{proof}
By Corollary~\ref{cor:corona},
\[
\gamma(G\circ H)=\gamma(G)+n\bigl(\gamma(H)+\delta(H)\bigr),\qquad
\gamma(H\circ G)=\gamma(H)+\mu\bigl(\gamma(G)+\delta(G)\bigr).
\]
The two quantities are equal if and only if $\gamma(G)-\mu\gamma(G)+n\gamma(H)-\gamma(H)=\mu\delta(G)-n\delta(H)$, which is the stated equation.
\end{proof}

\begin{corollary}\label{cor:sym}
Let $G$ and $H$ be connected graphs of orders $n$ and $\mu$.
\begin{enumerate}
\item[(a)] If $\delta(G)=\delta(H)=0$, then $\gamma(G\circ H)=\gamma(H\circ G)$ if and only if $(n-1)\gamma(H)=(\mu-1)\gamma(G)$.
\item[(b)] If $n=\mu\ge2$ and $\delta(G)=\delta(H)$, then $\gamma(G\circ H)=\gamma(H\circ G)$ if and only if $\gamma(G)=\gamma(H)$.
\end{enumerate}
\end{corollary}

\begin{proof}
Part (a) is Theorem~\ref{thm:sym} with $\delta(G)=\delta(H)=0$. For (b), the equation of Theorem~\ref{thm:sym} becomes $(n-1)(\gamma(H)-\gamma(G))=n(\delta(G)-\delta(H))=0$, and $n-1>0$.
\end{proof}

\begin{example}\label{ex:k6k33}
Let $G=K_6$ and $H=K_{3,3}$. Then $\gamma(G)=\gamma(H)=1$, and $K_1+K_6=K_7$ and $K_1+K_{3,3}=K_{3,3,1}$ both have genus $1$, so $\delta(G)=\delta(H)=0$. Both graphs have order $6$, so Corollary~\ref{cor:sym}(b) gives $\gamma(G\circ H)=\gamma(H\circ G)$, and indeed both equal $1+6\cdot1=7$.
\end{example}

Coronas of graphs of different genus can also have equal genus, and there are infinitely many such pairs. To construct them we use pendant vertices.

\begin{lemma}\label{lem:pendant}
Let $G$ be a connected graph, let $x\in V(G)$, and let $G'$ be obtained from $G$ by adding a new vertex $u$ and the edge $ux$. Then $\gamma(G')=\gamma(G)$ and $\gamma(K_1+G')=\gamma(K_1+G)$.
\end{lemma}

\begin{proof}
The graph $G'$ is the union of $G$ and the edge $ux$, which share only the vertex $x$, so $\gamma(G')=\gamma(G)+\gamma(K_2)=\gamma(G)$ by Theorem~\ref{thm:BHKY}. Let $\apex$ be the dominating vertex of $K_1+G$. In $K_1+G'$ the vertex $u$ is adjacent exactly to $x$ and $\apex$, and $\apex x$ is an edge of $K_1+G$. Take an embedding of $K_1+G$ in $S_g$ with $g=\gamma(K_1+G)$ and choose a face incident with the edge $\apex x$. Place $u$ in this face and join it to $\apex$ and to $x$ inside the face; this is possible because both lie on the boundary of the face. This embeds $K_1+G'$ in $S_g$, so $\gamma(K_1+G')\le\gamma(K_1+G)$. The reverse inequality follows from Lemma~\ref{lem:sub}.
\end{proof}

\begin{proposition}\label{prop:family}
For $n\ge5$ let $G_n$ be the graph obtained from $K_5$ by attaching $n-5$ pendant vertices to one vertex, and let $H_n$ be obtained from $K_4$ by attaching $n-3$ pendant vertices to one vertex. Then $|V(G_n)|=n$, $|V(H_n)|=n+1$, $\gamma(G_n)=1\ne0=\gamma(H_n)$, and
\[
\gamma(G_n\circ H_n)=\gamma(H_n\circ G_n)=n+1.
\]
\end{proposition}

\begin{proof}
The orders are clear. Applying Lemma~\ref{lem:pendant} once for each pendant vertex, $\gamma(G_n)=\gamma(K_5)=1$, $\gamma(K_1+G_n)=\gamma(K_6)=1$, $\gamma(H_n)=\gamma(K_4)=0$ and $\gamma(K_1+H_n)=\gamma(K_5)=1$. By Corollary~\ref{cor:corona},
\[
\gamma(G_n\circ H_n)=1+n\cdot1=n+1,\qquad \gamma(H_n\circ G_n)=0+(n+1)\cdot1=n+1. \qedhere
\]
\end{proof}

In this family $\delta(G_n)=0$ and $\delta(H_n)=1$, and both sides of the equation of Theorem~\ref{thm:sym} equal $-n$.

The corona is also not associative as an operation on graphs: if $|V(H)|=h$ and $|V(K)|=k\ge1$, then $|V((G\circ H)\circ K)|=n(1+h)(1+k)$ and $|V(G\circ(H\circ K))|=n(1+h(1+k))$, and these differ by $nk$. The genera can differ as well.

\begin{example}\label{ex:assoc}
Let $G=K_5$, $H=K_1$ and $K=K_4$. Since $K_1+K_1=K_2$ is planar, $\gamma(G\circ H)=\gamma(K_5)=1$, and $G\circ H$ has $10$ vertices. Because $K_1+K_4=K_5$,
\[
\gamma\bigl((G\circ H)\circ K\bigr)=1+10\,\gamma(K_5)=11.
\]
On the other hand $H\circ K=K_1\circ K_4=K_5$, so $G\circ(H\circ K)=K_5\circ K_5$ and
\[
\gamma\bigl(G\circ(H\circ K)\bigr)=\gamma(K_5)+5\,\gamma(K_6)=1+5=6.
\]
\end{example}

\subsection{A model with cut vertices, and graphs containing a generalized corona}
A backbone graph $G$ with a cluster $H_i$ attached to each vertex $v_i$, every device of $H_i$ being linked to $v_i$, is a generalized corona. By Theorem~\ref{thm:main} the genus splits into a term for the backbone and one independent term for each cluster. Table~\ref{tab:cluster} gives the cost $\gamma(K_1+H)$ of one cluster for common types.

\begin{table}[ht]
\centering
\caption{The value of $\gamma(K_1+H)$ for common graphs $H$.}\label{tab:cluster}
\begin{tabular}{@{}lc@{\qquad}lc@{}}
\toprule
$H$ & $\gamma(K_1+H)$ & $H$ & $\gamma(K_1+H)$\\
\midrule
tree, path or cycle & $0$ & $K_{3,3}$ & $1$\\
$K_4,\ K_5,\ K_6$ & $1$ & $K_{4,4}$ & $2$\\
$K_7$ & $2$ & wheel, prism, antiprism & $1$\\
Petersen graph $\Pet$ & $2$ &  & \\
\bottomrule
\end{tabular}
\end{table}

\begin{example}\label{ex:network}
Let $G=C_6$ and let the six vertices carry the graphs $K_4$, $K_7$, $C_5$, $\Pet$, $K_{3,3}$ and $P_4$. By Theorem~\ref{thm:main} and Table~\ref{tab:cluster},
\[
\gamma\bigl(C_6\circ(K_4,K_7,C_5,\Pet,K_{3,3},P_4)\bigr)=0+1+2+0+2+1+0=6.
\]
The graph has $6+4+7+5+10+6+4=42$ vertices and $101$ edges. The cycle and the outerplanar graphs $C_5$ and $P_4$ contribute nothing; all six handles come from the four dense graphs, each contributing independently.
\end{example}

The graph can be drawn without crossings on the sphere if and only if $G$ is planar and every $H_i$ is outerplanar (Corollary~\ref{cor:outer}). Monotonicity gives lower bounds for graphs that merely contain a generalized corona.

\begin{proposition}\label{prop:packing}
Let $X$ be a connected graph, let $G$ be a connected subgraph of $X$ with vertices $v_1,\dots,v_n$, and let $A_1,\dots,A_n$ be pairwise disjoint non-empty subsets of $V(X)\setminus V(G)$ with $A_i\subseteq N_X(v_i)$. Then
\[
\gamma(X)\ \ge\ \gamma(G)+\sum_{i=1}^n\gamma\bigl(K_1+X[A_i]\bigr).
\]
If $X$ is obtained from a generalized corona $Y$ by adding $q$ edges between vertices of $Y$, then $\gamma(Y)\le\gamma(X)\le\gamma(Y)+q$.
\end{proposition}

\begin{proof}
Put $H_i=X[A_i]$. The generalized corona $G\circ(H_1,\dots,H_n)$ is a connected subgraph of $X$: its edges are edges of $G$, edges inside the sets $A_i$, and the edges $v_ix$ with $x\in A_i$, all of which lie in $X$. The first inequality follows from Lemma~\ref{lem:sub} and Theorem~\ref{thm:main}. The second statement follows from Lemma~\ref{lem:sub} and from Lemma~\ref{lem:edge} applied $q$ times.
\end{proof}

For instance, if $q$ edges between vertices of different clusters are added to the graph of Example~\ref{ex:network}, the genus of the result lies between $6$ and $6+q$.

\section{Open problems}\label{sec:problems}

\begin{problem}
Determine $\delta(H)$ for further infinite families of non-planar graphs. The generalized Petersen graphs and the hypercubes $Q_n$ are natural candidates, and Table~\ref{tab:computed} gives data. By Proposition~\ref{prop:deltaK}, $\delta$ is known for $K_p$ and $K_{p,p}$.
\end{problem}

\begin{problem}
Is there a better bound than $c-1$ for planar $H$ in Proposition~\ref{prop:cover}, in terms of the vertex-face structure of $H$? The dodecahedron shows that $c-1$ is not sharp. More generally, how does $\gamma(K_1+H)$ depend on the faces needed to cover the vertices in embeddings of $H$ of higher genus?
\end{problem}

\begin{problem}
White's conjecture states that $$\gamma(K_{l,m,n})=\ceil{(l-2)(m+n-2)/4}$$ for $l\ge m\ge n$. A proof for all $K_{p,q,1}$ would determine $\gamma(G\circ K_{p,q})$ completely (Proposition~\ref{prop:unbal}).
\end{problem}

\appendix
\section{Rotation systems}\label{app:rot}

Each table gives, for every vertex, the cyclic order of its neighbours. Faces are obtained by the tracing rule of Section~\ref{sec:prelim}.

\begin{table}[ht]
\centering
\caption{A rotation system for the Petersen graph in the torus (Lemma~\ref{lem:petersen}).}\label{tab:rotP}
\begin{tabular}{@{}cl@{\qquad}cl@{}}
\toprule
vertex & rotation & vertex & rotation\\
\midrule
$0$ & $(1,4,5)$ & $5$ & $(0,7,8)$\\
$1$ & $(0,2,6)$ & $6$ & $(1,9,8)$\\
$2$ & $(1,3,7)$ & $7$ & $(2,5,9)$\\
$3$ & $(2,4,8)$ & $8$ & $(3,6,5)$\\
$4$ & $(0,3,9)$ & $9$ & $(4,6,7)$\\
\bottomrule
\end{tabular}
\end{table}

\begin{table}[ht]
\centering
\caption{A rotation system for $K_1+\Pet$ in the double torus (Theorem~\ref{thm:petersen}); $\apex$ is the added vertex.}\label{tab:rotKP}
\begin{tabular}{@{}cl@{\qquad}cl@{}}
\toprule
vertex & rotation & vertex & rotation\\
\midrule
$0$ & $(\apex,1,4,5)$ & $6$ & $(1,9,\apex,8)$\\
$1$ & $(0,\apex,2,6)$ & $7$ & $(\apex,2,5,9)$\\
$2$ & $(1,\apex,3,7)$ & $8$ & $(3,6,\apex,5)$\\
$3$ & $(2,\apex,4,8)$ & $9$ & $(4,\apex,6,7)$\\
$4$ & $(0,3,\apex,9)$ & $\apex$ & $(5,7,8,6,9,4,3,2,1,0)$\\
$5$ & $(\apex,0,7,8)$ & & \\
\bottomrule
\end{tabular}
\end{table}

\begin{table}[ht]
\centering
\caption{The twelve faces of the embedding in Table~\ref{tab:rotKP}, as closed walks.}\label{tab:facesKP}
\begin{tabular}{@{}cl@{}}
\toprule
length & boundary walks\\
\midrule
$3$ & $(0,\apex,5)$, $(0,1,\apex)$, $(1,2,\apex)$, $(2,3,\apex)$, $(3,4,\apex)$, $(4,9,\apex)$, $(6,\apex,9)$, $(6,8,\apex)$\\
$5$ & $(0,5,7,9,4)$, $(2,7,5,8,3)$\\
$6$ & $(0,4,3,8,6,1)$\\
$10$ & $(1,6,9,7,\apex,8,5,\apex,7,2)$\\
\bottomrule
\end{tabular}
\end{table}



\begin{thebibliography}{99}

\bibitem{Alpert73}
S.~R. Alpert, \emph{The genera of amalgamations of graphs}, Trans.\ Amer.\ Math.\ Soc.\ \textbf{178} (1973), 1--39.

\bibitem{BHKY62}
J.~Battle, F.~Harary, Y.~Kodama, and J.~W.~T. Youngs, \emph{Additivity of the genus of a graph}, Bull.\ Amer.\ Math.\ Soc.\ \textbf{68} (1962), 565--568.

\bibitem{BCHK16}
S.~Beyer, M.~Chimani, I.~Hedtke, and M.~Kotrb\v{c}\'ik, \emph{A practical method for the minimum genus of a graph: models and experiments}, Experimental Algorithms (SEA 2016), Lecture Notes in Computer Science, vol.~9685, Springer, 2016, pp.~75--88.

\bibitem{Brinkmann22}
G.~Brinkmann, \emph{A practical algorithm for the computation of the genus}, Ars Math.\ Contemp.\ \textbf{22} (2022), no.~4, \#P4.01.

\bibitem{ChartrandHarary67}
G.~Chartrand and F.~Harary, \emph{Planar permutation graphs}, Ann.\ Inst.\ H.\ Poincar\'e Sect.~B (N.S.) \textbf{3} (1967), 433--438.

\bibitem{Craft98}
D.~L. Craft, \emph{On the genus of joins and compositions of graphs}, Discrete Math.\ \textbf{178} (1998), 25--50.

\bibitem{DGH81}
R.~W. Decker, H.~H. Glover, and J.~P. Huneke, \emph{The genus of the $2$-amalgamations of graphs}, J.\ Graph Theory \textbf{5} (1981), 95--102.

\bibitem{ES07}
M.~N. Ellingham and C.~Stephens, \emph{The nonorientable genus of joins of complete graphs with large edgeless graphs}, J.\ Combin.\ Theory Ser.~B \textbf{97} (2007), 827--845.

\bibitem{ES09}
M.~N. Ellingham and C.~Stephens, \emph{The orientable genus of some joins of complete graphs with large edgeless graphs}, Discrete Math.\ \textbf{309} (2009), 1190--1198.

\bibitem{ESZ06}
M.~N. Ellingham, C.~Stephens, and X.~Zha, \emph{The nonorientable genus of complete tripartite graphs}, J.\ Combin.\ Theory Ser.~B \textbf{96} (2006), 529--559.

\bibitem{FruchtHarary70}
R.~Frucht and F.~Harary, \emph{On the corona of two graphs}, Aequationes Math.\ \textbf{4} (1970), 322--325.

\bibitem{GrossTucker87}
J.~L. Gross and T.~W. Tucker, \emph{Topological Graph Theory}, Wiley-Interscience, New York, 1987; Dover reprint, 2001.



\bibitem{KMR08}
K.~Kawarabayashi, B.~Mohar, and B.~Reed, \emph{A simpler linear time algorithm for embedding graphs into an arbitrary surface and the genus of graphs of bounded tree-width}, Proc.\ 49th Annual IEEE Symposium on Foundations of Computer Science, 2008, pp.~771--780.

\bibitem{Kurauskas17}
V.~Kurauskas, \emph{On the genus of the complete tripartite graph $K_{n,n,1}$}, Discrete Math.\ \textbf{340} (2017), 508--515.


\bibitem{LvChen19}
S.~Lv and Y.~Chen, \emph{Constructing a minimum genus embedding of the complete tripartite graph $K_{n,n,1}$ for odd $n$}, Discrete Math.\ \textbf{342} (2019), 3017--3024.

\bibitem{MaRen19}
D.~Ma and H.~Ren, \emph{The orientable genus of the join of a cycle and a complete graph}, Ars Math.\ Contemp.\ \textbf{17} (2019), 223--253.

\bibitem{Miller87}
G.~L. Miller, \emph{An additivity theorem for the genus of a graph}, J.\ Combin.\ Theory Ser.~B \textbf{43} (1987), 25--47.

\bibitem{Mohar99}
B.~Mohar, \emph{A linear time algorithm for embedding graphs in an arbitrary surface}, SIAM J.\ Discrete Math.\ \textbf{12} (1999), 6--26.

\bibitem{MoharThomassen01}
B.~Mohar and C.~Thomassen, \emph{Graphs on Surfaces}, Johns Hopkins University Press, Baltimore, 2001.

\bibitem{Laali16}
A.~R.~F. Laali, H.~H. Seyyedavadi, and D.~Kiani, \emph{Spectra of generalized corona of graphs}, Linear Algebra Appl.\ \textbf{493} (2016), 411--425.


\bibitem{Ringel65}
G.~Ringel, \emph{Das Geschlecht des vollst\"andigen paaren Graphen}, Abh.\ Math.\ Sem.\ Univ.\ Hamburg \textbf{28} (1965), 139--150.

\bibitem{Ringel74}
G.~Ringel, \emph{Map Color Theorem}, Grundlehren der mathematischen Wissenschaften, vol.~209, Springer, Berlin, 1974.

\bibitem{RY68}
G.~Ringel and J.~W.~T. Youngs, \emph{Solution of the Heawood map-coloring problem}, Proc.\ Nat.\ Acad.\ Sci.\ U.S.A.\ \textbf{60} (1968), 438--445.

\bibitem{Sharma17}
R.~Sharma, B.~Adhikari, and A.~Mishra, \emph{Structural and spectral properties of corona graphs}, Discrete Appl.\ Math.\ \textbf{228} (2017), 14--31.

\bibitem{StahlBeineke77}
S.~Stahl and L.~W. Beineke, \emph{Blocks and the nonorientable genus of graphs}, J.\ Graph Theory \textbf{1} (1977), 75--78.

\bibitem{StahlWhite76}
S.~Stahl and A.~T. White, \emph{Genus embeddings for some complete tripartite graphs}, Discrete Math.\ \textbf{14} (1976), 279--296.

\bibitem{Thomassen89}
C.~Thomassen, \emph{The graph genus problem is NP-complete}, J.\ Algorithms \textbf{10} (1989), 568--576.

\bibitem{White69}
A.~T. White, \emph{The genus of the complete tripartite graph $K_{mn,n,n}$}, J.\ Combin.\ Theory \textbf{7} (1969), 283--285.

\bibitem{Youngs63}
J.~W.~T. Youngs, \emph{Minimal imbeddings and the genus of a graph}, J.\ Math.\ Mech.\ \textbf{12} (1963), 303--315.

\end{thebibliography}
\end{document}